\documentclass[11pt]{article}
\usepackage{amsfonts}
\usepackage{amscd}
\usepackage{amsmath}
\usepackage{amsthm}

\newtheorem{thm}{Theorem}
\newtheorem{prop}[thm]{Proposition}

\newcommand \al{\alpha}

\newcommand\de{\delta}

\newcommand\et{\eta}

\newcommand\ro{\rho}

\newcommand\ph{\varphi}

\newcommand\om{\omega}

\newcommand\Om{\Omega}

\def\ZZ{\mathbb Z}

\def\RR{\mathbb R}
\def\CC{\mathbb C}

\def\TT{\mathbb T}
\newcommand\oo{{\infty}}

\renewcommand\o{\circ}
\renewcommand\div{\on{div}}
\newcommand\x{\times}
\newcommand\on{\operatorname}

\newcommand\ad{\on{ad}}

\newcommand\grad{\on{grad}}

\newcommand\Diff{\on{Diff}}

\newcommand\id{\on{id}}
\newcommand\Der{\on{Der}}

\newcommand\g{\mathfrak g}

\newcommand\h{\mathfrak h}

\newcommand\e{\mathfrak e}
\newcommand\X{\mathfrak X}

\date{ }
\begin{document}

\title{Geodesics on extensions of Lie groups and stability;
the superconductivity equation}
\author{Cornelia Vizman \\\it \small West University of Timisoara,
Department of Mathematics\\ 
\it\small Bd. V.Parvan 4, 1900-Timisoara, Romania\\
\it \small e-mail: vizman@math.uvt.ro}
\maketitle

\begin{abstract}
The equations of motion of an ideal charged fluid,
respectively the superconductivity equation
(both in a given magnetic field) are showed to be geodesic equations 
of a general, respectively a central extension of the group
of volume preserving diffeomorphisms with right invariant metrics.
For this, quantization of the magnetic flux is required. 
We do curvature computations in both cases in order to get 
informations about the stability.
\end{abstract}




\section{Introduction}

Important partial differential equations were obtained 
as geodesic equations on diffeomorphism groups with right invariant 
$L^2$ or $H^1$ metrics: 
Euler equation of motion of an incompressible ideal fluid
is the geodesic equation on the volume preserving diffeomorphism
group with $L^2$ metric \cite{A} \cite{MEF},
Burger's equation is the geodesic equation on $\Diff(S^1)$
with   $L^2$ metric \cite{AK}, the Korteweg-de-Vries equation 
is the geodesic equation on 
the Virasoro-Bott group with $L^2$ metric \cite{OK}, 
the Camassa-Holm shallow water equation \cite{CH} is the geodesic equation on
$\Diff(S^1)$ and on the Virasoro-Bott group with  
$H^1$ metrics \cite{M2} \cite{K} \cite{CM}, 
the averaged Euler equation is the geodesic
equation on 
the volume preserving diffeomorphism group with $H^1$ metric
\cite{MRS} \cite{S}, the equations of ideal magneto-hydrodynamics are geodesic
equations on the semidirect product of the group of volume preserving
diffeomorphisms and the linear space of divergence free vector fields 
with $L^2$ metric \cite{MRW} \cite{ZK}.

The Lagrangian (exponential) instability of geodesics is determined by 
the negativity of the sectional curvature.
In this way Arnold \cite{A} showed the instability in most directions 
of Euler equation for ideal flow, Shkoller \cite{S} showed that 
averaged Euler equation is more stable than Euler equation for ideal flow,
Misiolek \cite{M2} \cite{M1} has results on the stability of the 
Korteweg-de-Vries and Camassa-Holm equations and Zeitlin and Kambe
\cite{ZK} for ideal hydrodynamics. 

In this paper we are concerned with the equations of motion 
of an ideal charged fluid in a given magnetic field
$u_t=-\nabla_u u-\rho u\times B-\grad p$, 
$\rho_t=-d\rho.u$, $\div u=0$ and the superconductivity equation
(obtained for $\rho=1$) as geodesic equations.

We consider the central extension by the 1-torus $\TT$, due to Ismagilov
\cite{I1}\cite{I2},
of the group of exact volume preserving diffeomorphisms 
of a compact Riemannian manifold $M$,
corresponding to the Lichnerowicz Lie algebra cocycle
$\om(X,Y)=\int_M\et(X,Y)\mu$, where $\et$ is a closed 2-form on $M$. 
This extension exists only if 
the cohomology class of $\et$ is integral.
In the 3-dimensional case $\et$ corresponds to the magnetic
vector field $B$ and the integrability condition on $\et$
is equivalent to flux quantization of magnetic field:
$\int_S(B\cdot n)dS\in\ZZ$ for every closed surface $S$ in $M$. 
In the case $H^{n-1}(M,\RR)=0$, the group of exact volume
preserving diffeomorphisms coincides with the group of volume
preserving diffeomorphisms, and the superconductivity equation is proved 
to be the geodesic equation on the central extension with right invariant
$L^2$-metric.

To drop out the condition $H^{n-1}(M,\RR)=0$ we can consider
a non-central extension of the group of volume preserving diffeomorphisms:
the group of volume preserving automorphisms of a principal
$\TT$-bundle $P$ over $M$ possessing a principal connection with curvature
form $\et$. Again we have to impose the integrality condition on $[\et]$.
We consider the right invariant $L^2$-metric coming from the
Kaluza-Klein metric on $P$.
Then the geodesic equations are the equations of motion of
an ideal charged fluid.
For $\rho=1$ we get again the superconductivity equation as 
a geodesic equation.

To get information on the 
Lagrangian stability of the superconductivity equation, we compute
the sectional curvature in the case
of the flat 3-torus in both settings (general and central extension)
and the result is: 
the superconductivity equation is more stable than 
Euler equation for ideal flow.  

Like in the book of Arnold and Khesin \cite{AK}, 
the approach followed throughout the paper is formal. For example a proof
of the existence of geodesic flow on Fr\'echet manifolds is not 
available in the literature: for the intricacies of the rigorous geometric 
approach to this field see Hamilton's paper about the inverse function
theorem of Nash and Moser \cite{H}.
The ILH (inverse limit of Hilbert) Lie group setting and well posedness
of the Cauchy problem for the superconductivity equation
will be treated elsewhere.

\section{Superconductivity equation}

Let $M$ be an $n$-dimensional 
compact Riemannian manifold with Levi Civita connection $\nabla$
and volume form $\mu$. Let $B$ be an $(n-2)$ vector field on $M$ 
(i.e. $B\in C^\oo(\wedge^{n-2}TM)$) 
such that $\et=(-1)^{n-2}i_B\mu$ is a closed
two-form. The cross product of a vector field $X$ with $B$ is 
the vector field
$X\x B=(i_{X\wedge B}\mu)^\sharp=(i_X\et)^\sharp$.

Then we can write the generalized
equations of motion of an ideal charged fluid 
for a time dependent divergence free vector field $u$ 
and a time dependent real smooth function $\rho$ on $M$:
\begin{gather}
u_t=-\nabla_uu-\rho u\times B-\grad p\nonumber\\
\rho_t=-d\rho.u,
\end{gather}
where the index $t$ denotes partial derivation by $t$.
When at $t=0$ the function $\rho$ is a constant, say $\rho=1$,
then $\rho$ is a constant at every $t$ and
we get the generalized superconductivity equation
\begin{equation}
u_t=-\nabla_uu-u\x B-\grad p.
\end{equation}

In case $M$ is a 3-dimensional manifold, $B$ is simply a vector field and 
the condition $\et=-i_B\mu$ closed is equivalent to the condition $B$
divergence free. Moreover, the cross product is the usual vector product.
Then equation (1) models the motion of an ideal charged fluid
in a given magnetic field $B$: $u$ represents the velocity field 
and $\rho$ the charge density.
The superconductivity equation (2) models the  
motion of a high density electronic gas in a magnetic field: 
$u$ represents the velocity field of 
the electronic gas.
Its similarity to Euler equation 
comes from the fact that at high density, because of the repelling
of particles, an electron gas behaves like a fluid. 

Zeitlin \cite{Z} and Roger \cite{R} already showed that
the superconductivity equation
can be regarded as the Euler equation on the central extension of 
the Lie algebra of divergence free
vector fields on $M$ by the Lichnerowicz cocycle.
Poisson structures for superconductors are given in \cite{HK}.

If the closed 2-form $\et$ is integral we will obtain both 
equations (1) and (2) as geodesic equations on extensions
of the group of volume preserving diffeomorphisms.
The integrality condition on $\et$ translates into a quantization of 
the magnetic flux: the flux of $B$ across any closed surface is integral
\begin{equation}
\int_S(B\cdot n)dS\in\ZZ.
\end{equation}
This magnetic field seems to be produced by Dirac monopols.

\section{Right invariant metrics on Lie groups}

In this paragraph we give expressions for the geodesic equation, Levi-Civita 
covariant derivative and curvature for Lie groups with right invariant 
metrics (see \cite{MR} for a nice presentation of this subject).

Let $G$ be a Lie group with Lie algebra $\g$. Let $\ro_x$ be
the right translation by $x$. Any right invariant bounded Riemannian 
metric on $G$ is determined by its value at the identity $<,>:\g\x\g\to\RR$, 
a positive definite bounded inner product on $\g$.
Let $g:I\to G$ be a smooth curve and $u:I\to\g$ its right logarithmic 
derivative (the velocity field in the right trivialization) 
$u(t)=T\ro_{g(t)^{-1}}.g'(t)$.
In terms of $u$ the geodesic equation for $g$ has the expression
\begin{equation}
u_t=-\ad(u)^\top u,
\end{equation}
where $\ad(X)^\top$ is the adjoint of $\ad(X)$ with respect to $<,>$, 
if this adjoint does exist. 

The right trivialization induces an isomorphism 
$R:C^\oo(G,\g)\to\mathfrak{X}(G)$ given by $R_X(x)=T\ro_x.X(x)$. In terms of 
this isomorphism, the Levi-Civita covariant derivative is
\begin{equation}
\nabla^G_X Y=dY.R_X+\frac12\ad(X)^\top Y+\frac12\ad(Y)^\top X-\frac12\ad(X)Y,
\end{equation}
for $X,Y\in C^\oo(G,\g)$.

The sectional curvature ${\cal K}$ and the Riemannian curvature ${\cal R}$ 
are related by 
\begin{equation}
{\cal K}(X,Y)=\frac{<{\cal R}(X,Y)Y,X>}{<X,X><Y,Y>-<X,Y>^2},
\end{equation} 
so the sign of the expression 
$<{\cal R}(X,Y)Y,X>$ determines the sign of the sectional curvature. 
Its expression in the right trivialization is
\begin{align}
<{\cal R}(X,Y)Y,X>&=\frac14\Vert\ad(X)^\top Y+\ad(Y)^\top X\Vert^2
-\frac34\Vert\ad(X)Y\Vert^2\nonumber\\
&-<\ad(X)^\top X,\ad(Y)^\top Y>\\
&-\frac12<\ad(X)^\top Y,\ad(X)Y>
-\frac12<\ad(Y)^\top X,\ad(Y)X>.\nonumber
\end{align}

{\bf An example}: \cite{A} \cite{MEF}
Let $G=\Diff_{vol}(M)$ be the group of volume preserving diffeomorphisms 
of a compact Riemannian manifold $(M,g)$ with induced volume form $\mu$,
Levi Civita covariant 
derivative $\nabla$ and Riemannian curvature tensor $R$.
Let $\g=\mathfrak{X}_{vol}(M)$ 
the Lie algebra of divergence free vector fields.
We consider the right invariant metric on $G$ given by the $L^2$ inner 
product $<X,Y>=\int_M g(X,Y)\mu$. The transpose of $\ad(X)$ 
is $\ad(X)^\top Y=P(\nabla_X Y+(\nabla X)^\top Y)$,
the geodesic equation in terms of the right logarithmic derivative
is Euler equation for ideal flow
\begin{equation}
u_t=-\nabla_uu-\grad p,\quad\div u=0,
\end{equation}
the covariant derivative for right invariant vector fields is 
$\nabla^G_XY=P\nabla_XY$ and the curvature 
\begin{align}
<{\cal R}(X,Y)Y,X>
&=<R(X,Y)Y,X>\\
&+<Q\nabla_XX,Q\nabla_YY>-\Vert Q\nabla_XY\Vert^2,\nonumber
\end{align}
with $P$ and $Q$ the orthogonal projections 
on the spaces of divergence free respectively gradient vector fields.

Let $\TT^2$ be the two-dimensional torus and $X_k$ the Hamiltonian
vector field with Hamiltonian (stream function) $\cos kx$, $k\in\ZZ^2$, 
$x\in\TT^2$. 
The curvature of $\Diff_{vol}(\TT^2)$ in any two-dimensional direction 
containing the direction $X_k$ is non-positive.

\section{Geodesics on general Lie group extensions}

Consider an exact sequence of homomorphisms of Lie groups 
\begin{equation}
1\to H\to E\to G\to 1
\end {equation}
and the corresponding exact sequence of homomorphisms of Lie algebras
\begin{equation}
0\to\h\to\e\to\g\to 0,
\end{equation}
i.e. the Lie group $E$ is an extension of $G$ by $H$. 
A section $s:\g\to\h$ induces the following mappings:
\begin{align}
&b:\g\to\Der(\h), \quad b(X)Y=[s(X),Y]\\
&\om:\g\x\g\to\h, \quad\om(X_1,X_2)=[s(X_1),s(X_2)]-s([X_1,X_2])
\end{align}
with properties
\begin{align}
&[b(X_1),b(X_2)]-b([X_1,X_2])=\ad(\om(X_1,X_2))\\
&\sum_{cycl}\om([X_1,X_2],X_3)=\sum_{cycl}b(X_1)\om(X_2,X_3).
\end{align}
The Lie algebra structure on the extension $\e$,
identified as vector space with $\g\oplus\h$ via the section $s$, 
in terms of $b$ and $\om$ is \cite{AMR}
\begin{equation}
[(X_1,Y_1),(X_2,Y_2)]=([X_1,X_2],[Y_1,Y_2]+b(X_1)Y_2-b(X_2)Y_1+\om(X_1,X_2)
\end{equation}
Special cases: when $\om=0$ we get semidirect products, when $\h$ is
abelian and $b=0$ we get central extensions with Lie algebra cocycle $\om$.

We consider also the right invariant metric on the Lie group $E$ given at 
the identity by the positive inner product on $\e$:
\begin{equation}
<(X_1,Y_1),(X_2,Y_2)>_\e=<X_1,X_2>_\g+<Y_1,Y_2>_\h,
\end{equation}
where $<,>_\g$ and $<,>_\h$ are positive definite inner products on 
$\g$ and $\h$ such that the transposes $\ad(X)^\top:\g\to\g$ and
$\ad(Y)^\top:\h\to\h$ exist for any $X\in\g$ and $Y\in\h$. 
To write the 
transpose of $\ad(X,Y)$ in the extended Lie algebra, 
we have to impose further conditions.
We suppose the transpose
$b(X)^\top:\h\to\h$ exists for any $X\in\g$ and
there exist maps $h:\h\to L(\g)$ linear 
(actually $h$ takes values
in the space of skew-adjoint operators on $\g$) and $l:\h\x\h\to\g$ bilinear,
defined by the relations
\begin{gather}
<\om(X_1,X_2),Y>_\h=<h(Y)X_1,X_2>_\g\\
<b(X)Y_1,Y_2>_\h=<l(Y_1,Y_2),X>_\g.
\end{gather}
Then the transpose of $\ad(X,Y)$ is
\begin{align}
\ad(X_1,Y_1)^\top(X_2,Y_2)=&(\ad(X_1)^\top X_2+h(Y_2)X_1-l(Y_1,Y_2),\nonumber\\
&\ad(Y_1)^\top Y_2+b(X_1)^\top Y_2)
\end{align}

\begin{prop}
With the conditions above, the geodesic equation on the Lie group extension
$E$ with right invariant metric, written in terms of the right 
logarithmic derivative $(u,\rho)$ with $u:I\to\g$ and $\rho:I\to\h$ is:
\begin{gather}
u_t=-\ad(u)^\top u-h(\rho)u+l(\rho,\rho)\nonumber\\
\rho_t=-\ad(\rho)^\top\rho-b(u)^\top\rho
\end{gather}
\end{prop}

We are interested in the special case when
$b(X)$ is skew-adjoint for all $X\in\g$: $<b(X)Y_1,Y_2>+<Y_1,b(X)Y_2>=0$.
This will be called the isometric case, like in \cite{Vi}, 
where the semidirect product with this special condition is studied.

\begin{prop}
On an isometric extension of $\g$ by abelian $\h$, 
the geodesic equation in terms of the right logarithmic derivative $(u,\rho)$
is 
\begin{gather}
u_t=-\ad(u)^\top u-h(\rho)u\nonumber\\
\rho_t=b(u)\rho
\end{gather}
and the curvature tensor at the identity for the extended group is:
\begin{gather}
<{\cal R}^E((X_1,Y_1),(X_2,Y_2))(X_2,Y_2),(X_1,Y_1)>
=<{\cal R}^G(X_1,X_2)X_2,X_1>\nonumber\\
-<b(X_1)Y_2-b(X_2)Y_1,\om(X_1,X_2)>
-\frac34\Vert\om(X_1,X_2)\Vert^2\nonumber\\
+\frac14\Vert h(Y_2)X_1+h(Y_1)X_2\Vert^2-<h(Y_1)X_1,h(Y_2)X_2>\\
+<\om(X_1,\nabla^G_{X_1}X_2),Y_2>+<\om(X_2,\nabla^G_{X_2}X_1),Y_1>\nonumber\\
-<\om(X_2,\nabla^G_{X_1}X_1),Y_2>-<\om(X_1,\nabla^G_{X_2}X_2),Y_1>,\nonumber
\end{gather}
where ${\cal R}^G$ denotes the Riemannian curvature at the identity in $G$ and 
${\cal R}^E$ in the extended group $E$.
\end{prop}

In the particular case of an isometric semidirect product of $\g$ 
with an abelian Lie 
algebra $\h$, the only non-zero term at the right hand side is
$<{\cal R}^G(X_1,X_2)X_2,X_1>$.

{\bf Geodesics and curvature on one-dimensional extensions of Lie groups}

Next we consider another special case: one-dimensional central extensions 
$\tilde G$ of $G$ with Lie 
algebra $\tilde g$ defined by the Lie algebra cocycle
$\om:\g\x\g\to\RR$. 
Then the Lie algebra bracket on $\tilde\g$ is
$[(X_1,a_1),(X_2,a_2)]=([X_1,X_2],\om(X_1,X_2))$.
On the extended group we consider the right invariant metric 
given at the identity by the following 
positive definite inner product on $\tilde\g$:
$<(X_1,a_1),(X_2,a_2)>=<X_1,X_2>+a_1a_2$.

Suppose the transpose $\ad(X)^\top:\g\to\g$ with respect to 
the inner product exists for all $X\in\g$ and the Lie algebra cocycle 
can be written in the form $\om(X_1,X_2)=<k(X_1),X_2>$ 
with $k:\g\to\g$; it follows $k$ is skew-adjoint operator and $k(X)=h(1)X$.

\begin{prop}
With the conditions above, the geodesic equation 
on the one-dimensional central extension $\tilde G$,
in terms of the right logarithmic derivative $(u,a)$ is
\begin{equation}
u_t=-\ad(u)^\top u-ak(u), a\in\RR.
\end{equation}
The sign of the sectional curvature in the extended group is given
by the sign of
\begin{gather}
<\widetilde{\cal R}((X_1,a_1),(X_2,a_2))(X_2,a_2),(X_1,a_1)>=
<{\cal R}(X_1,X_2)X_2,X_1>\nonumber\\
-\frac34\om(X_1,X_2)^2
+\frac14\Vert k(a_1X_2-a_2X_1)\Vert^2\\
-\om(\nabla^G_{a_1X_2-a_2X_1}X_1,X_2)+\om(\nabla^G_{a_1X_2-a_2X_1}X_2,X_1),
\nonumber
\end{gather}
where ${\cal R}$ denotes the Riemannian curvature at the identity in $G$ and 
$\widetilde{\cal R}$ in the extended group $\tilde G$.
\end{prop}

With this formula applied to the Virasoro-Bott 
extension, considering $L^2$ and $H^1$ metrics,
we can recover the results of Misiolek on the sign of the sectional curvature
and so on the 
Lagrangian (exponential) stability of the Korteweg-de-Vries and 
Camassa-Holm equations \cite{M2} \cite{M1}.

\section{Lichnerowicz cocycle}

Let $M$ be a compact Riemannian manifold with Levi Civita connection $\nabla$
and volume form $\mu$. 
Any closed two-form $\et$ on $M$ determines a Lie algebra cocycle, called 
Lichnerowicz cocycle, on the Lie algebra of divergence free vector fields
$\mathfrak{X}_{vol}(M)$:
\begin{equation}
\om(X,Y)=\int_M\et(X,Y)\mu
\end{equation}

Does it exist a corresponding Lie group cocycle on the group of volume
preserving diffeomorphisms $\Diff_{vol}(M)$? The existence of a smooth
local group cocycle $c$ integrating $\om$ follows from \cite{N}, where this
result is proved in general for Lie groups modeled over sequentially 
complete locally convex spaces.

When the cohomology class of $\et$ is integral, 
results of Ismagilov \cite{I1} \cite{I2}
imply the existence of a global group extension 
of the subgroup of exact volume
preserving diffeomorphisms $\Diff_{exact}(M)$ by the torus $\TT$,
corresponding 
to a Lichnerowicz Lie algebra extension determined by $\et$.
A volume preserving 
diffeomorphism is called exact if there is an isotopy $\ph_t\in\Diff_{vol}(M)$
from $\ph_0=\id_M$ to $\ph_1=\ph$ such that  
$\frac{d}{dt}\ph_t=X_t\o\ph_t$ and $i(X_t)\mu$ is an exact $(n-1)$-form.
In case $H^{n-1}(M,\RR)=0$, the subgroup of exact diffeomorphisms coincides
with the identity component of $\Diff_{vol}(M)$.
 
\begin{prop}
Let $M$ be a compact manifold with $H^{n-1}(M,\RR)=0$ and 
$[\et]\in H^2(M,\ZZ)$. Then the generalized 
superconductivity equation 
\begin{equation}
u_t=-\nabla_uu-u\times B-\grad p.
\end{equation}
is a geodesic
equation on a central extension of the group of volume
preserving diffeomorphisms with right invariant $L^2$ metric.
\end{prop}

\begin{proof}
The transpose of $\ad(X)$ with respect to
the $L^2$ inner product on the Lie algebra of divergence free vector fields
exists: $\ad(X)^\top X=P\nabla_XX$.
and the Lichnerowicz cocycle is of the required form  
$\om(X,Y)=<k(X),Y>$ with $k(X)=P(X\times B)$.
Then proposition 3 can be applied with initial condition $a=1$.
\end{proof}

\section{Volume preserving automorphisms of\\
a principal $\TT$-bundle}

Let $M$ be an $n$-dimensional Riemannian manifold with volume form $\mu$ and 
$\et$ a closed differential two-form with integral cohomology class
$[\et]$. Then there exists a principal $\TT$-bundle $\pi:P\to M$
and a principal connection form $\al\in\Om^1(P)$ having curvature $\et$.
The associated Kaluza-Klein metric on $P$, defined at a point $x\in P$
by
\begin{equation}
\kappa_x(\tilde X,\tilde Y)=g_{\pi(x)}
(T_x\pi.\tilde X,T_x\pi.\tilde Y)+\al_x(\tilde X)\al_x(\tilde Y),
\end{equation}
determines the volume form $\tilde\mu=\pi^*\mu\wedge\al$ on $P$.

The group of volume preserving automorphisms of the
principal bundle $P$, denoted by
$E=\Diff_{vol}(P)^\TT=\Diff_{vol}(P)\cap\Diff(P)^\TT$
is an extension
of the group of volume preserving diffeomorphisms $G=\Diff_{vol}(M)$
by the abelian gauge group $H=C^\oo(M,\TT)$. 
The corresponding exact sequence of Lie
algebras is
\begin{equation}
0\to C^\oo(M,\RR)\to\X_{vol}(P)^\TT\to\X_{vol}(M)\to 0.
\end{equation}
The principal connection $\al$ provides us with a section $s(X)=X^H$,
$X^H$ representing the horizontal lift of the divergence free vector field
$X$ on $M$. It is easy to verify that the horizontal lift $X^H$
is again divergence free, with respect to the volume form 
$\tilde\mu$ on $P$. 
The identification of $\X_{vol}(P)^\TT$ with $\X_{vol}(M)\oplus C^\oo(M,\RR)$
(as vector spaces only) via the section $s$ is $(X,f)\mapsto X^H+(f\o\pi)\xi$,
where $\xi$ denotes the fundamental vector field of the principal torus action
(in particular $\al(\xi)=1$).

The mapping $b$ induced by $s$ is minus the usual
action of vector fields as derivations on the algebra of smooth functions
$b(X)f=-df.X$
and $\om$ is the restriction of the curvature two-form $\et$  
to divergence free vector fields, since the curvature measures the deviation of
the horizontal bundle from integrability:
$(\et(X_1,X_2)\o\pi)\xi=[X_1^H,X_2^H]-[X_1,X_2]^H$. 

The $L^2$ inner product on $\e=\X_{vol}(P)^\TT$ coming from the 
Kaluza-Klein metric on $P$ is easily shown to be 
the sum of the $L^2$ inner products 
on $\g=\X_{vol}(M)$ and $\h=C^\oo(M,\RR)$:
\begin{equation}
<(X_1,f_1),(X_2,f_2)>_\e=\int_M g(X_1,X_2)\mu+\int_M f_1f_2\mu
\end{equation}
(here we impose the condition that the integral of the 1-form $\al$
along the fibers of $P$ is 1; otherwise we get a constant factor
at the right side of the equality above).

\begin{prop}
The geodesic equation on the group of divergence free automorphisms 
of the principal $\TT$-bundle $P$ with right invariant $L^2$ metric 
coming from the Kaluza-Klein metric on $P$, written in terms
of the right logarithmic derivative $(u,\rho)$,
$u:I\to\X_{vol}(M)$, $\rho:I\to C^\oo(M,\RR)$,
is:
\begin{gather}
u_t=-\nabla_uu-\rho u\x B-\grad p\nonumber\\
\rho_t=-d\rho.u,
\end{gather}
the equations of motion of an ideal charged fluid.
\end{prop}

\begin{proof}
All the conditions required in proposition 2
are fulfilled by this extension. The mapping
$h:C^\oo(M,\RR)\to L(\X_{vol}(M))$ is $h(f)X=P(fX\x B)$
and the mapping $l:C^\oo(M,\RR)\x C^\oo(M,\RR)\to\X_{vol}(M)$
is $l(f_1,f_2)=P(f_1\grad f_2)$. We have again denoted by $P$ the 
orthogonal projection on the space of divergence free
vector fields on $M$. The transpose $\ad(X)^\top X=P(\nabla_XX)$
and $b(X)$ is skew-adjoint. Now  
proposition 2 can be applied.
\end{proof}

The geodesic equation on $\Diff_{vol}(P)^\TT$,
for the initial condition $(u_0,a)$ with 
$u_0\in\X_{vol}(M)$ and $a\in\RR$,
is $u_t=-\nabla_uu-au\x B-\grad p$, like for the central
extension corresponding to the Lichnerowicz cocycle. 
So, if $a=1$, we obtain the superconductivity equation as a 
geodesic equation also 
for arbitrary compact Riemannian manifolds $M$ (i.e. without
the condition $H^{n-1}(M,\RR)=0$ we had to impose in section 5), 
especially in the interesting case of the 
flat 3-torus $\TT^3$. 
We remark that the vector subspace $\{(X,a):X\in\X_{vol}(M),a\in\RR\}$
is not a Lie subalgebra of $\X_{vol}(P)^\TT$.

\section {Curvature computations}

The expressions in propositions 2 and 3 determine the 
sign of the sectional curvature
and we get informations on the Lagrangian stability of geodesics. 
A geodesic is Lagrangian stable if all geodesics with sufficiently 
close initial conditions at time zero remain close for all $t\ge 0$ \cite{M3}.
This is a notion of stability referring to the motion, 
not to the velocity field. 

The sign of the sectional curvatures of the same two-dimensional plane
spanned by $(X_1,a_1)$, $(X_2,a_2)$ with $a_1,a_2\in\RR$ 
can differ for the two extensions 
of $\Diff_{vol}(M)$: the central extension $\tilde G$ and the general
extension $E=\Diff_{vol}(P)^{\TT}$. 
Using $h(a)X=ak(X)$ for $a\in\RR$ and $X\in\X_{vol}(M)$ we get
\begin{gather}
<{\cal R}^E((X_1,a_1),(X_2,a_2))(X_2,a_2),(X_1,a_1)>=
<{\cal R}^G(X_1,X_2)X_2,X_1>\nonumber\\
-\frac34\Vert\et(X_1,X_2)\Vert^2
+\frac14\Vert k(a_1X_2-a_2X_1)\Vert^2\\
-\om(\nabla^G_{a_1X_2-a_2X_1}X_1,X_2)+\om(\nabla^G_{a_1X_2-a_2X_1}X_2,X_1),
\nonumber
\end{gather}
where $\nabla^G_X Y=P\nabla_XY$ and the formula for ${\cal R}^G$ is (9).
This expression is very 
similar to the expression of the curvature in the central 
extension, only instead of $\om(X_1,X_2)^2=(\int_M\et(X_1,X_2)\mu)^2$
we get in this case $\Vert\et(X_1,X_2)\Vert^2=\int_M\et(X_1,X_2)^2\mu$.

{\bf The flat 3-torus}

Next we do curvature computations for $M=\TT^3$, the 3-torus,
in both settings: central extension by $\RR$ and general extension
by $C^\oo(M,\RR)$,
although in the first case only the existence of a local
group cocycle integrating the Lichnerowicz cocycle is known
($H^2(\TT^3)\ne 0$).

Let $M=\TT^3=\RR^3/(2\pi\ZZ)^3$ be the flat 3-torus.
The Fourier basis is $e_k(x)=e^{ik\cdot x}$, where $k\in\ZZ^3$. We complexify
the Lie algebra of divergence free vector fields on the torus, as well as 
the inner product, commutator and Levi Civita connection. A vector field
on the torus is written in the Fourier basis as $u=\sum_k u_ke_k$ with
$u_k\in\CC^3$. The reality condition is $u_{-k}=\bar u_k$ and the divergence
free condition is $k\cdot u_k=0$. Moreover
\begin{align}
[u,v]&=i((l\cdot u_k)v_l-(k\cdot u_l)v_k)e_{k+l}\\
Q(u)&=\frac1{\vert k\vert^2}(k\cdot u_k)ke_k\\
\nabla_uv&=i(l\cdot u_k)v_le_{k+l}.
\end{align}

A 3-dimensional version of Arnold's non-positivity result for the curvature
of the group of volume preserving diffeomorphisms on the 2-torus 
was obtained in \cite{NHK}: for 
$U_p=u_pe_p+u_{-p}e_{-p}$ and $X=\sum_k v_ke_k$,
\begin{equation}
<{\cal R}(U_p,X)X,U_p>=
-\Vert Q\nabla_{U_p}X\Vert^2\le 0.
\end{equation}

\begin{prop}
Let $B=B_0$ be a constant vector field on the 3-torus.
Let $U_p=u_pe_p+u_{-p}e_{-p}$ and $V_p=v_pe_p+v_{-p}e_{-p}$ be divergence 
free vector fields on the 3-torus with the additional
condition $u_p$ and $v_p$ real and orthogonal. 
Then the sectional curvature 
of the two-dimensional plane spanned by 
$(U_p,1)$ and $(V_p,0)$
on the central extension $\tilde G$ 
(respectively on the general extension $E$) of $\Diff_{vol}(M)$ 
is positive if and only if 
$\vert u_p\vert^2<\frac1{6(2\pi)^3}$ (respectively $\vert u_p\vert^2<\frac13$).
\end{prop}

\begin{proof}
We observe that
$\nabla_{U_p}V_p=0$, so
in the curvature formula (14) 
applied to $(U_p,1),(V_p,0)$, only 2 
terms are non-zero:
\begin{equation}
<\tilde{\cal R}((U_p,1),(V_p,0))(V_p,0),(U_p,1)>
=-\frac34\om(U_p,V_p)^2+\frac14\Vert k(V_p)\Vert^2.
\end{equation}
Using the Fourier coefficients, we get the following expressions 
of the cocycle $\om$ and of the mapping $k$:
\begin{gather}
\om(u,v)^2=(\int_{\TT^3}\mu(B,u,v)\mu)^2=(2\pi)^6(\sum\mu(B_k,u_l,v_m)
\de(k+l+m))^2\\
\Vert\et(u,v)\Vert^2=\int_{\TT^3}\mu(B,u,v)^2\mu=(2\pi)^3\sum
\mu(B_k,u_l,v_m)^2\de(k+l+m)
\end{gather}
and
\begin{equation}
k(u)=P(u\times B)=((u_l\x B_k)-\frac1{\vert k+l\vert^2}
((k+l)\cdot(u_l \x B_k))(k+l))e_{k+l}.
\end{equation}
From the hypothesis, $(p,u_p,v_p)$ is an orthogonal basis of $\RR^3$. In 
the associated orthonormal basis the coordinates of $B_0$ are 
denoted by $(a_1,a_2,a_3)$.
Let $\tilde{\cal R}$ (respectively ${\cal R}^E$) be the curvature 
tensor of the central (respectively general) extension of the group
of volume preserving diffeomorphisms. Then
\begin{align}
<\tilde{\cal R}&((U_p,1),(V_p,0))(V_p,0),(U_p,1)>
=-\frac34(2\pi)^6(\mu(B_0,u_p,\bar v_p)\nonumber\\
&+\mu(B_0,\bar u_p,v_p))^2
+\frac12(2\pi)^3(B_0\x v_p-\frac1{\vert p\vert^2}\mu(p,B_0,v_p)p)^2\\
&=\frac12(2\pi)^3a_1^2\vert v_p\vert^2(1-6(2\pi)^3\vert u_p\vert^2).\nonumber
\end{align}
So the sectional curvature on the central extension is positive if and only if
$\vert u_p\vert^2<\frac1{6(2\pi)^3}$.
An analogous computation for ${\cal R}^E$ gives the condition 
$\vert u_p\vert^2<\frac13$.
\end{proof}

In conclusion, the superconductivity equation on the 3-torus 
is more stable 
than Euler equation for ideal flow. 

{\it Acknowledgments}: 
I address special thanks to the referee for suggesting
to consider the $L^2$-metric coming from the
Kaluza-Klein metric associated with the principal connection.
Part of this paper was done during my visit at the Institute
of Mathematics in Erlangen;
I gratefully acknowledge Prof. K. Strambach and the staff of the Institute 
for their hospitality.


\end{document}